\documentclass[a4paper, 10pt, conference]{ieeeconf} 
\pdfoutput=1
\IEEEoverridecommandlockouts
\usepackage{amsfonts,bm}
\usepackage{amsmath}
\usepackage{amssymb}
\usepackage[pdftex]{graphicx}
\DeclareGraphicsExtensions{.pdf}

\def\R{\mathbb{R}}
\def\C{\mathbb{C}}

\def\ci{\mathfrak{i}}
\def\ud{\textrm{d}}
\def\e{\textrm{e}}

\def\prt#1#2{\frac{\partial #1}{\partial #2}}

\newcommand{\be}{\begin{equation}}
\newcommand{\ee}{\end{equation}}
\newcommand{\bea}{\begin{eqnarray}}
\newcommand{\eea}{\end{eqnarray}}
\newcommand{\ben}{\begin{displaymath}}
\newcommand{\een}{\end{displaymath}}
\newcommand{\bean}{\begin{eqnarray*}}
\newcommand{\eean}{\end{eqnarray*}}

\newtheorem{remark}{Remark}

\title{\LARGE \bf
Open-Loop Control Design via Parametrization\\
Applied in a Two-Level Quantum System Model
}
\author{Markku Nihtil\"a,  {\it Senior Member, IEEE}
\thanks{M. Nihtil\"a is with Department of Physics and Mathematics,
        University of Eastern Finland, Kuopio Campus, 
        Yliopistonranta 1, P.O.Box 1627, FI-70211, Kuopio, Finland
        {\tt\small Markku.Nihtila@uef.fi}}
}
\begin{document}
\maketitle
\thispagestyle{empty}
\pagestyle{empty}
\begin{abstract}
In the design of quantum computing devices of the future the basic element 
is the qubit. It is a two-level quantum system which may describe population transfer 
from one steady-state to another controlled by a coherent laser field. A four-dimensional 
real-variable differential equation model is constructed from the complex-valued two-level 
model describing the wave function of the system. The state transition matrix of the model 
is constructed via the Wei-Norman technique and Lie algebraic methodology. The idea of 
parametrization using flatness-based control, is applied to construct feasible input--output 
pairs of the model. This input drives the state of the system from the given initial state to the 
given final state in a finite time producing the corresponding output of the pair. The population transfer is obtained by nullifying part of the state vector via careful selection of the parameter functions. A preliminary simulation study completes the paper.
\end{abstract}
\section{INTRODUCTION}
In quantum mechanical framework deterministic bits "1" and "0" are substituted by the {\it qubit} \cite{Wiki_Qubit}.
The qubit is a composition of the  pure states "1" and "0". This composition means that
the actual state of the qubit is not exactly "1" or "0" but a combination of these. In the measurement, 
however, the outcome is always one of the two possibilities "1" or "0". The qubit can be represented as a point 
on the surface of a sphere, so-called {\it Bloch sphere}, see e.g. \cite{Wiki_Qubit}. If one wants to save 
information into a qubit, then the key peoblem is to drive the qubit from one state to another. Then one 
arrives at the description of the qubit as a dynamic differential equation system, the controls of which are the parameters of the driving laser field.

In quantum computation the qubit forms a basic element for building up multi-qubit computing elements of future quantum computers, see \cite{Nature08a}. Then a key problem is to drive the qubit from one stable level to another.

Molecular excitation,  i.e. driving of an ensemble of molecules from one locally stable steady state to another is one alternative for a qubit structure. This type of systems are controlled by using coherent light. Based on laser technology shorter and shorter coherent pulses can be  
generated for controlling molecular excitation, see \cite{Nature08a}--\cite{Sonia09a}. The goal is to direct 
molecular reactions towards unprobable but desirable direction \cite{bandra02}-\cite{LeBris00}. Then 
nonlinear and more and more sophisticated control methods are needed 
for properly designing durations and forms of the control pulses. In classical N-level problems  the system to be controlled can be modelled by using ordinary  $2N$-dimensional differential equation systems. Due to 
femto- and picosecond scale pulses feedback is not in general applicable in 
the control design for these systems. Flatness-based control, see \cite{fliess90}-\cite{fliess99} \& \cite{Levine06}, is then an ideal methodology for 
open-loop design relevant in quantum control problems. 

Due to the fact that the bilinear quantum control systems are not controllable in the whole Euclidian space $\R^{2N}$ the methodology applied here is called {\it parametrization}.

This two-level quantum control problem and some related studies have been carried out by several authors earlier, too, see \cite{IEEE03}, \cite{IEEE04}, \& \cite{IEEE02}. Especially, in \cite{IEEE02} a very similar approach as ours is used.

However, we start from the basic definition of differential flatness.
The system 
\bea
\frac{dx}{dt}&=&f(x,u) ; \ \ x(t) \in  
{\mathbb R}^n , \ \ u(t) \in  \mathbb{R}^m  \label{state}
\eea

\noindent is called differentially flat if there exists algebraic 
functions (\cite{fliess90}) ${\mathcal A}$, ${\mathcal B}$, ${\mathcal C}$, 
and finite integers $\alpha$, $\beta$, and $\gamma$ 
such that for any pair $(x, u)$ of inputs and controls, satisfying 
the dynamics (\ref{state}), there exists a function $z$, 
called a flat (or linearizing) output,
such that the following equations are satisfied
\bea\label{ex,uu}
&&x(t)= {\mathcal A}(z, \dot z, \, \dots, z^{(\alpha)}) \nonumber \\
&&u(t)= {\mathcal B}(z, \dot z, \, \dots, z^{(\beta)}) \\
&&z(t)= {\mathcal C}(x, u,
 \dot u, \, \dots, u^{(\gamma)}). \nonumber
\eea
The actual output $y$, which is not present in the definition of flatness, may have the dependence
$$y(t) = h(x(t),u(t))$$
for some given output function $h$.
In parametrization procedure, due to uncontrollability, the last equation in (\ref{ex,uu}) for $z(t)$ is neither constructed nor applied.

From the standard finite-state Schr\"odinger equation of two energy levels a four-dimensional real-variable differential equation model is obtained. The Wei-Norman technique is used in the construction according to \cite{weinorman64}. The exponential representation of the transition matrix of the system includes three base functions, two of which serve as the parameter functions. In this framework the initial and final states can be defined corresponding to the two levels of the original system model. Then parametrization design is applied for explicitly calculating the parameter functions, which in turn give the desired input--output pairs.
\section{SYSTEM MODELS}
Population transfer in a two-level quantum system, see \cite{Boscain02}, can be described by the time-dependent Schr\"odinger equation, {\it i.e.} by the dynamics
\bea \label{mato}
\ci \frac{d \tilde \psi}{d t} =  \tilde H(t)\, \tilde \psi, \ \ \tilde H(t)=
\left[ \begin{array}{cc}
E_1 & \Omega(t) \\
\Omega^*(t) & E_2
\end{array} \right],
\eea
where the modified Planck's constant $\hbar=\frac{h}{2 \pi}$ has been scaled to $\hbar=1$, and $\ci =\sqrt{-1}.$ The wavefunction $ \tilde \psi : \R \to \C^2$ has the probabilistic interpretation, in the sense that 
\bea \label{proba1}
\Vert \tilde \psi (t) \Vert^2 = |\tilde \psi_1(t)|^2 + |\tilde \psi_2(t) |^2 = 1, \, \forall\, t \in \R,
\eea
where $\tilde \psi = (\tilde \psi_1, \tilde \psi_2)$. The control is given by  $\Omega : \R  \to \C$, and $\Omega^*$ is the complex conjugate of $\Omega$. $E_1$ and $E_2$ are the energy levels. The unitary transformation $\tilde \psi \mapsto \psi$ and $\Omega \mapsto u$ by
\bea
\tilde \psi (t)&=& U(t)\,  \psi (t), \\
U(t) &=&
\left[ \begin{array}{cc}
e^{-\ci E_1 t} & 0 \\
0 & e^{-\ci E_2 t}
\end{array} \right] \\[0.2cm]
u(t) &=& e^{-\ci(E_2- E_1) t}\, \Omega(t)
\eea
 transforms (\ref{mato}) to
\bea \label{ilman}
\ci \frac{d  \psi}{d t} &=&   H(t)\,  \psi, \\ 
H(t)&=&\left[ \begin{array}{cc}
0 & u(t) \\
u^*(t) & 0
\end{array} \right].
\eea
The componentwise representation
\bea
&&\psi (t)=\psi_1(t)  \left[ \begin{array}{c}
1 \\
0
\end{array} \right]
 + \psi_2 (t)   \left[ \begin{array}{c}
0 \\
1
\end{array} \right] 
\eea
converts (\ref{ilman}) to the dynamics
\bea \label{psii}
\begin{array}{l}
\dot \psi_1 = - \ci \, u \, \psi_2, \\
\dot \psi_2 = - \ci \, u^* \, \psi_1. 
\end{array}
\eea
By using the real-valued decompositions
\bea
\begin{cases}
\psi_1 = x_1+ \ci \, x_2 \\
\psi_2 = x_3+ \ci \, x_4 \\
\, \, u \, =\,  u_1 +\ci \, u_2
\end{cases}
\eea
one obtains a state-variable representation
\bea \label{tila4}
\left[ \begin{array}{l}
\dot x_1 \\
\dot x_2 \\
\dot x_3 \\
\dot x_4
\end{array} \right]
= \left[ \begin{array}{rr}
x_4 &   x_3 \\
- x_3 & x_4 \\
x_2 &  - x_1 \\
- x_1 & -x_2
\end{array} \right]
\left[ \begin{array}{l}
u_1 \\
u_2
\end{array} \right] 
\eea
or in another form
\bea \label{statex}
\dfrac{\ud x}{\ud  t} = \big(u_1 F_1  + u_2 F_2\big) x,
\eea
\bea
x = \left[ x_1 \,x_2 \,x_3 \,x_4\right]^T, 
\eea
\bea
F_1 &=&\left[ \hspace*{-0.1cm}\begin{array}{rrrr}
0 & 0 & 0 & 1 \\
0 & 0 & -1 & 0 \\
0 & 1 & 0 & 0 \\
-1 & 0 & 0 & 0
\end{array} \right], \\
F_2 &=& \left[\hspace*{-0.1cm}\begin{array}{rrrr}
0 & 0 & 1 & 0 \\
0 & 0 & 0 & 1 \\
-1 & 0 & 0 & 0 \\
0 & -1 & 0 & 0
\end{array} \right].
\eea
 The constraint (\ref{proba1}) is converted into the form
\bea
\sum_{k=1}^4 x_k^2 = 1.
\eea
\remark The matrices $F_1$ and $F_2$ together with their Lie product $2 F_3 =[F_1, F_2] =
 F_1F_2-F_2F_1$ form a Lie algebra. This can be used as a basis for differential geometric considerations of the control system (\ref{tila4}). However, the elementary approach applied in this paper is sufficient for our parametrization purposes.

\section{WEI-NORMAN REPRESENTATION}

The Lie algebra of the matrices $F_1, \, F_2$, and $F_3$ is three-dimensional with the relations
\bea
&&[ F_1, F_2 ] = 2 F_{3}\, ,  \\[0.1cm]
&&[ F_2, F_3 ] = 2 F_{1} \, , \\[0.1cm]
&&[ F_3, F_1 ] = 2 F_{2}.
\eea
\bea
F_3=\left[
\begin{array}{rrrr}
0 & -1 & 0 & 0 \\
1 & 0 & 0 & 0 \\
0 & 0 & 0 & 1 \\
0 & 0 & -1 & 0
\end{array} \right]
\eea
Due to the linear structure of the system model (\ref{statex}) with respect to the state $x$, the state transition matrix of the system, denoted by $\Phi$, and which relates the values of the state according to\\[-0.5cm]
\bea \label{Phi01}
x(t) = \Phi (t,0)\, x(0)
\eea
can be written as a product of exponentials
\bea \label{transmatrix}
\Phi (t,0)= 
 \e^{g_1 F_1}\,\e^{g_2 F_2}\,\e^{g_3 F_3}\, ,
\eea
where the exponentials are defined by the absolutely converging infinite series
\bea
\e^{g_i F_i} = \sum_{k=0}^{\infty} \dfrac{1}{k!}\, g_i^k F_i^k\, , \ \ i=1, 2, 3.
\eea
The state transition matrix satisafies the following initial-value problem (IVP1)\\[-0.5cm]
\bea \label{ivalprobl}
\prt{}{t}\Phi (t,0) &=& F(t)\, \Phi (t,0); \ \ \Phi (0,0) =I, \\
F(t) &=& u_1(t) F_1 + u_2(t) F_2 + 0 \cdot F_3. \label{fii01}
\eea
The technique we are using is nowadays called Wei-Norman technique according to the paper of Wei and Norman \cite{weinorman64}.
Substitution of the (\ref{transmatrix}) to the IVP1 gives
\bea \label{fii02}
\prt{}{t}\Phi \hspace*{-0.25cm}&=&\hspace*{-0.25cm}\dot g_1 F_1\, \Phi + \dot g_2\, \e^{g_1F_1} F_2\, \e^{-g_1F_1}\, \Phi  \\
\hspace*{-0.25cm}&+&\hspace*{-0.25cm}\dot g_3\, \e^{g_1F_1} \e^{g_2F_2}F_3\,  \e^{-g_2F_2} \e^{-g_1F_1} \, \Phi \, . \nonumber
\eea
By using (several times) the Campbell-Baker-Hausdorff formula for square matrices $A$ and $B$ of the same dimension\\[-0.5cm]
\bea
\e^A\, B \, \e^{-A}&=&B + [A, B] + [A, [A, B]]/2! \\
&+& [A, [A, [A, B]]]/3! + \cdots \nonumber
\eea
in the equation (\ref{fii02}) it can be represented in the form
\bea \label{CBH}
\prt{}{t}\Phi &=& \left[ f_1(t) F_1 +  f_2(t) F_2 +  f_3(t) F_3 \right]\, \Phi \\ [0.2cm]
f_1(t) &=&  \dot g_1 + \dot g_3 \sin(2g_2) \\[0.2cm]
f_2(t)&=& \dot g_2 \cos (2g_1) - \dot g_3 \cos (2g_2) \, \sin (2g_1) \\[0.2cm]
f_3(t)&=& \dot g_2 \sin (2g_1) + \dot g_3 \cos (2g_2)\, \cos (2g_1)
\eea
By comparing the coefficients of the $F_i$'s in (\ref{CBH}) and (\ref{ivalprobl})--(\ref{fii01}) one finally obtains a differential relation between the $g_i$'s and the controls $u_1$ and $u_2$ in the form of a matrix equation

\bea \label{WN01}
\left[ \begin{array}{ccc}
1 & 0 & \sin(2g_2)  \\[0.1cm]
0 & \cos (2g_1) & - \cos (2g_2) \, \sin (2g_1) \\[0.1cm]
0 & \sin (2g_1) &  \ \ \, \cos (2g_2)\, \cos (2g_1) 
\end{array} \right] \hspace*{-0.2cm}
\left[
\begin{array}{c}
\dot g_1 \\[0.1cm]
\dot g_2 \\[0.1cm]
\dot g_3 
\end{array} \hspace*{-0.10cm}\right]\hspace*{-0.18cm}=\hspace*{-0.18cm}
\left[
\begin{array}{c}
u_1 \\[0.1cm]
u_2 \\[0.1cm]
0
\end{array} \hspace*{-0.10cm}\right] \hspace*{-0.18cm} \\
{} \nonumber
\eea
the coefficient matrix being the same as in \cite{altafini02}, Eq. (3.7).
The relation $g \leftrightarrow u$ is invertible if the determinant of the coefficient matrix denoted by $\mathcal{D}$ is different from zero
\bea
|\mathcal{D}| = \cos( 2 g_2) \neq 0.
\eea
Then we have
\bea
\mathcal{D}^{-1} = \dfrac{1}{\cos( 2 g_2)} \times 
\eea
\vspace*{-0.2cm}
\bea
\left[ \begin{array}{ccc}
\cos (2g_2)  & - \sin (2g_1) \,\sin(2g_2)    & \cos (2g_1) \, \sin(2g_2)  \\
0 & \ \ \, \cos (2g_1) \, \cos (2g_2)  &  \sin (2g_1) \, \cos (2g_2)  \\
0 & - \sin (2g_1) &  \cos (2g_1) 
\end{array} \right] \nonumber 
\eea
\bea
g = \mathcal{D}^{-1} \tilde u \, ,
\eea
where $\tilde u$ and $g$ are defined by
\bea
\tilde u = \left[
\begin{array}{c}
u_1 \\[0.1cm]
u_2 \\[0.1cm]
0
\end{array} \right]\, , \ \ g =
\left[
\begin{array}{c}
\dot g_1 \\[0.1cm]
\dot g_2 \\[0.1cm]
\dot g_3 
\end{array} \right].
\eea

\section{MODEL PARAMETRIZATION}

Because the system has two (scalar) controls we can choose two of the three base functions $g_i$ freely corresponding to free selection of the two controls. The third base function has to be determined from the last equation of (\ref{WN01}). Parametrization actually means that the input--output pairs can be determined from the parameter functions without explicitly solving of the system equations according to Fig. 1.
\begin{figure}
\begin{center}
\includegraphics[width=8.0cm]{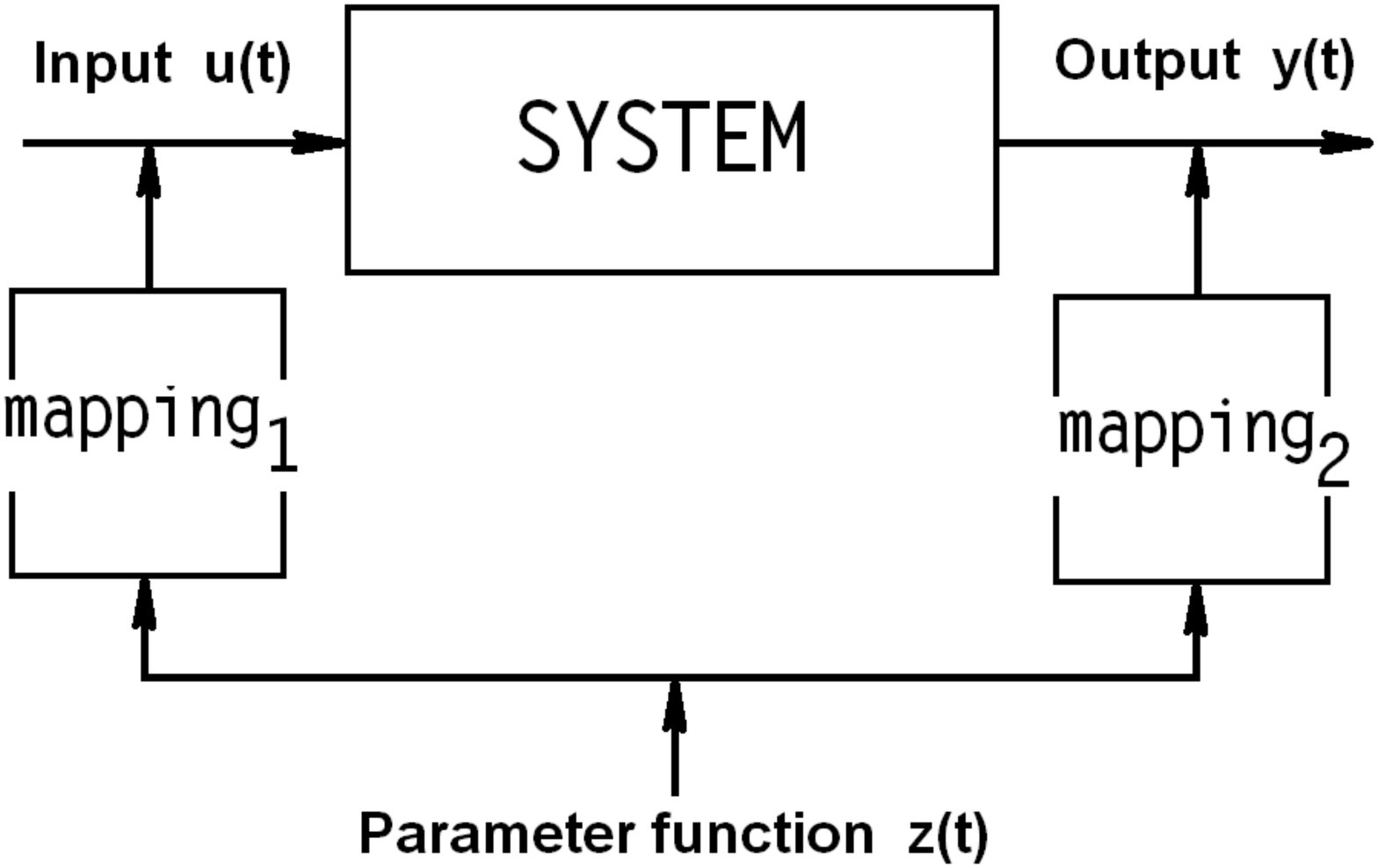}
\vspace*{-0.3cm}
\caption{Mappings 1 \& 2 give the corresponding input-output pairs $(u, y)$
without explicitly solving the system equations.}
\end{center}
\end{figure}
Due to the flatness-based design idea, computation of the third base function as well as of the controls must not include integrations as given by the equations (\ref{ex,uu}). Only differentiations are allowed. Consequently, based on the third equation in (\ref{WN01}), the base functions $g_2$ and $g_3$ are chosen as parameter functions. Then these are also so-called flat outputs, see \cite{fliess95}, denoted by $z=(z_1, z_2)=(g_2, g_3)$. The parametrization obtained in this way for $g_1$ and the controls
are given by
\bea
g_1 &=& \frac{1}{2}\,\arctan \Big[-\cos (2 g_2)\, \dfrac{\dot g_3}{\dot g_2}\, \Big] \label{derieq} \\[0.2cm]
u_1 &=& \dot g_1 + \dot g_3 \sin (2 g_2) \label{u1} \\[0.2cm]
u_2 &=& \sqrt{\dot g_2^2 + \dot g_3^2 \cos^2(2 g_2)}\, . \label{u2}
\eea
The state variables are calculated by using the state transition matrix equations (\ref{Phi01}) and  (\ref{transmatrix})
\ben
x(t) = \Phi (t, 0)x(0) =  \e^{g_1 F_1}\,\e^{g_2 F_2}\,\e^{g_3 F_3}\, x(0).
\een
\section{Control objective}
In population transfer problems from the level 1 corresponding to the situation 
\bea
|\psi_1(0)|^2=x_1(0)^2+x_2(0)^2=0
\eea
 to the level 2, where 
\bea
|\psi_2(T)|=x_3(T)^2+x_4(T)^2=0 \, ,
\eea
where $T$ is the transfer time, we can parametrize the partial trajectory by using a sufficiently smooth, but otherwise arbitrarily chosen, parametrization $x_1, x_2$ with the boundary conditions
\bea
x_1(0)^2+x_2(0)^2=0, \label{alku} \\
x_1(T)^2+x_2(T)^2=1. \label{loppu}
\eea
By dividing the state vector into two parts
\bea
x(t) &=& (\, w(t), v(t)\, ) \\
w(t) &=& (\, x_1(t), x_2(t)\,) \\
v(t) &=& (\, x_3(t), x_4(t)\,) 
\eea
we can represent the task of driving the state from the initial one to the final one in a finite time $T$ as follows
\bea
x(0)=\left[ \begin{array}{c}
0 \\
0 \\
x_{30} \\
x_{40}
\end{array} \right]  &\rightarrow &
\left[ \begin{array}{c}
x_{1T} \\
x_{2T}  \\
0 \\
0
\end{array} \right] = x(T)\\[0.2cm]
\left[ \begin{array}{c}
0 \\
0 \\
\sin \alpha \\
\cos \alpha
\end{array} \right]  &\rightarrow  &
\left[ \begin{array}{c}
\cos \beta \\
\sin \beta \\
0 \\
0
\end{array} \right] \hspace*{2.0cm}
\eea

We have chosen a specific parametrization for the initial and final values of the state, because the sum of the squares of the nonzero state components must be equal to 1 at the both ends of the planned trajectory.

\section{PARAMETRIZATION DESIGN}

The state transition equation $x(T)= \Phi (T, 0) x(0)$ can now be written in the form
\bea
\left[ \begin{array}{c}
w_T \\
0 \end{array} \right] 
= \left[ \begin{array}{cc}
A & B \\
C & D \end{array} \right] \left[ \begin{array}{c}
0 \\
v_0 \end{array} \right] 
\eea
\bea
\therefore \ \  \left\{
\begin{array}{r}
w_T = B v_0 \\
0 = D v_0.
\end{array}.  \right.
\eea
where $A, \, B, \, C,$ and $D$ are $2\hspace*{-0.06cm}\times\hspace*{-0.06cm}2$-blocks of the $4\hspace*{-0.06cm}\times\hspace*{-0.06cm}4$-dimensional state transition matrix $\Phi(T, 0)$.

For the state transition matrix
\bea \label{transmatrix02}
\Phi (t,0)= 
 \e^{g_1 F_1}\,\e^{g_2 F_2}\,\e^{g_3 F_3}
\eea
where the exponentials are defined by the series
\bea
\e^{g_i F_i} = \sum_{k=0}^{\infty} \dfrac{1}{k!}\, g_i^k F_i^k\, , i=1, 2, 3
\eea
we obtain the series representations in closed form
\bea
\e^{g_i F_i} =\cos g_i \, I + \sin g_i \, F_i
\eea
due to the fact that $F_i^2 = - I, \, i=1, 2, 3$, where $I$ is $4\hspace*{-0.06cm}\times\hspace*{-0.06cm}4$ identity matrix. Then the product of the three exponent functions is of the form
\bea
&&\Phi = (c_1 I + s_1 F_1) (c_2 I + s_2 F_2) (c_3 I + s_3 F_3) \\[0.2cm]
&&c_i = \cos g_i, \ \ s_i = \sin g_i, \ \ i=1, 2, 3.
\eea
Now the $D$-part and $B$-part of the transfer matrix $\Phi$ are given by
\bea
D &=& \left[ \begin{array}{cc}
d_1 & d_2 \\
d_3 & d_4 \end{array} \right] \\[0.1cm]
&=& c_1c_2 \left[ \begin{array}{rr}
c_3 & s_3 \\
-s_3 & c_3 \end{array} \right] - 
s_1s_2 \left[ \begin{array}{rr}
s_3 & -c_3 \\
c_3 & s_3 \end{array} \right] , \label{ehto1}\\[0.1cm]
B &=& \left[ \begin{array}{cc}
b_1 & b_2 \\
b_3 & b_4 \end{array} \right] \\[0.1cm]
&=& c_1s_2 \left[ \begin{array}{rr}
c_3 & s_3 \\
-s_3 & c_3 \end{array} \right] - 
s_1c_2 \left[ \begin{array}{rr}
s_3 & -c_3 \\
c_3 & s_3 \end{array} \right]. 
\eea
We must have $D =0$
due to the requirement
\mbox{$D v_0 =0$}
for arbitrary \mbox{$v_0 = (\, x_{30}, x_{40}\, ) $} satisfying the requirement \mbox{$x_{30}^2+ x_{40}^2 = 1$}. Then we have two alternatives in (\ref{ehto1}):
\bea
\begin{cases}
a) & c_1 =s_2 =0 \\
& \textrm{or} \\
b) & s_1=c_2 =0
\end{cases} \ \ \Rightarrow \ \ D =0 \ \ \therefore \ \ Dv_0=0.
\eea
These conditions are obtained from the two basic alternatives
\bea
a) \ \
\begin{cases}
\cos g_1(T)=0\, , & g_1(T)= \dfrac{\pi}{2} \\
\sin g_2(T)=0\, , & g_2(T)=  0,
\end{cases} 
\eea
\bea
b) \ \
\begin{cases}
\sin g_1(T)=0\, , & g_1(T)= 0 \\
\cos g_2(T)=0\, , & g_2(T)=  \dfrac{\pi}{2}.
\end{cases} 
\eea
In the case of the first alternative {\it a)} we have
\bea
\begin{cases}
s_1=  \sin g_1(T) = 1 \\
c_2 = \cos g_2(T) = 1 .
\end{cases}
\eea
Consequently,
\bea
B &=&- 
s_1c_2 \left[ \begin{array}{rr}
s_3 & -c_3 \\
c_3 & s_3 \end{array} \right],
\eea
\bea
\hspace*{-0.2cm}w_T &=& B v_0 = - \left[ \begin{array}{rr}
\sin g_3 & -\cos g_3 \\
\cos g_3 & \sin g_3 \end{array} \right] 
\left[ \begin{array}{c}
\sin \alpha  \\
\cos \alpha  \end{array} \right] \\[0.1cm]
&=& - \left[ \begin{array}{c}
\sin g_3 \sin \alpha  -\cos g_3 \cos \alpha \\
\cos g_3 \sin \alpha + \sin g_3 \cos \alpha \end{array} \right] 
\\[0.1cm]
&=& \left[ \begin{array}{c}
\cos (- g_3 - \alpha ) \\
 \sin (- g_3 - \alpha ) \end{array} \right] 
= \left[ \begin{array}{c}
\cos \beta \\
\sin \beta  \end{array} \right] .
\eea
\bea
\therefore \ \ g_3(T)= -\, ( \alpha + \beta ) .
\eea
In the same way the alternative {\it b)} can be solved giving
\bea
g_3(T)= \dfrac{\pi}{2} - ( \alpha + \beta ) .
\eea 
Due to trigonometric functions in the equations there are also other possibilities for the final values of $g_2$ and $g_3$ deviating by the multiples of $\pi$ or $2 \pi$. These possibilities need further considerations and are not studied here. We choose the alternative {\it b)} for the basis of our control design. So, we have to find sufficiently differentiable parameter 
functions $g_2$ and $g_3$, which together with the dependent basis function $g_1$
have to satisfy the boundary conditions
\bea \label{bcond}
\begin{cases}
g_1(0) = 0 \, , & g_1(T) = 0\, ; \\[0.2cm]
g_2(0) = 0  \, , & g_2(T) =\dfrac{\pi}{2}\, ; \\[0.2cm]
g_3(0) = 0   \, , & g_3(T) =\dfrac{\pi}{2}\, - (\alpha + \beta ) \, .
\end{cases} 
\eea
 The final value of $g_1$ depends on the derivatives of $g_2$ and $g_3$. This means that we have to adjust these derivatives via the equation ({\ref{derieq}) to agree with the requirement $g_1(T)= 0$.

Carefully planned and realized simulations are needed to confirm the feasibility of our parametrization approach.
\section{SIMULATION STUDY}
First preliminary simulation results demonstrate that the methodology developed actually drives the state of the system from the given initial state (level 1) to the given final state (level 2). A minimal parametrization for the parameter functions $g_2$ and $g_3$ were chosen without any specific optimization procedure. The only requirements are that the given boundary conditions (\ref{bcond}) are satisfied, and that the equation which gives the base function $g_1$ also gives the correct initial and final values for $g_1$. The following values were used in the simulations
\bea
\begin{cases}
\alpha &= -2 \pi/3 \, , \\
\beta&= \, \pi/3 \, , \\
T &= 10 \, .
\end{cases}
\eea
Then the final value for $g_3$ becomes
\bea
g_3(T)=\dfrac{5 \pi}{6} \approx 2.62 \, .
\eea
Because $g_2$ has to change from $0$ to $\pi/2$, we chose the linear function
\bea \label{g2}
g_2(t) = \dfrac{\pi }{2}\, \dfrac{t}{T}\approx 0.157\, t \, .
\eea
The boundary values
\bea
\begin{cases}
g_1(0)&=0 \\
g_1(T)&=0
\end{cases}
\eea
 are obtained when we choose 
\bea
\begin{cases}
\dot g_3(0)&=0 \\
\dot g_3(T)&=0
\end{cases}.
\eea
Then the third order polynomial suffices
\bea
g_3(t) = a_0 + a_1 \, \dfrac{t}{T}  +a_2 \Big(\dfrac{t}{T} \Big)^2 + a_3 \Big(\dfrac{t}{T} \Big)^3 \, .
\eea
The coefficients are obtained from the boundary conditions, giving finally
\bea \label{g3}
g_3(t) = \gamma \left\{3 \Big(\dfrac{t}{T} \Big)^2 - 2 \Big(\dfrac{t}{T} \Big)^3 \right\} \, , \ \ 
\gamma = \dfrac{\pi}{2}\, - (\alpha + \beta ) \, .
\eea
The binding condition
\bea
 \dot g_2 \sin (2g_1) + \dot g_3 \cos (2g_2)\, \cos (2g_1)=0
\eea
gives the base function $g_1$ for the given parameter functions (\ref{g2}) and (\ref{g3}).
The functions are depicted in Fig. 2 and 3. The controls were calculated by using the formulas (\ref{u1}) and (\ref{u2}). They are depicted in Fig. 4.
\begin{figure}
\begin{center}
\includegraphics[width=4.2cm]{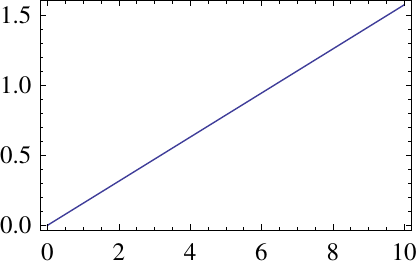}
\includegraphics[width=4.2cm]{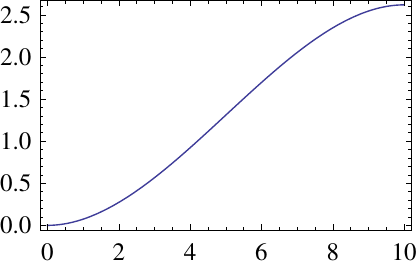}
\vspace*{-0.3cm}
\caption{The independent parameter functions: $g_2$ -- left, $g_3$ -- right.}
\end{center}
\vspace*{0.4cm}
\begin{center}
\includegraphics[width=4.2cm]{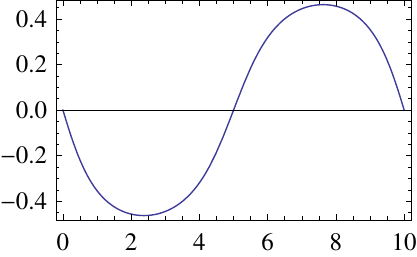}
\vspace*{-0.3cm}
\caption{The dependent base function $g_1$.}
\end{center}
\end{figure}
\vspace*{0.3cm}
\begin{figure}
\includegraphics[width=4.2cm]{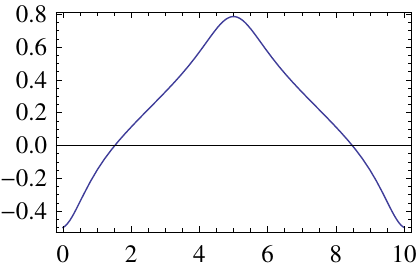}
\includegraphics[width=4.2cm]{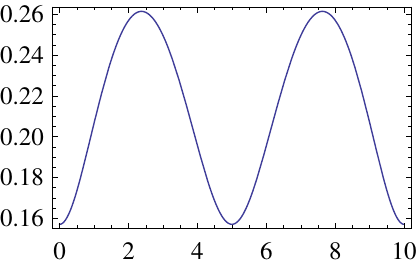}
\vspace*{-0.3cm}
\caption{The control variables $u_1(t)$ and $u_2(t)$ for $t \in [0, 10]$.}
\vspace*{0.3cm}
\begin{center}
\includegraphics[width=4.2cm]{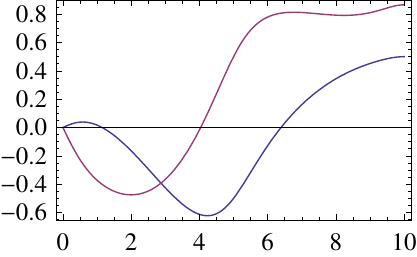}
\includegraphics[width=4.2cm]{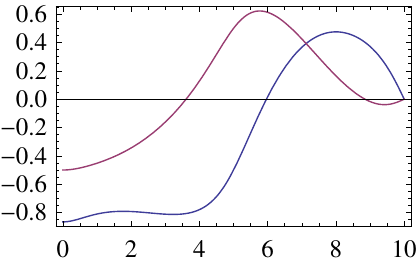}
\vspace*{-0.3cm}
\caption{The state variables obtained by simulating the system equations by using the given control functions. 
Left: $x_1$ -- lower, $x_2$ -- upper. Right: $x_3$ -- lower, $x_4$ -- upper.}
\end{center}
\vspace*{0.3cm}
\begin{center}
\includegraphics[width=4.2cm]{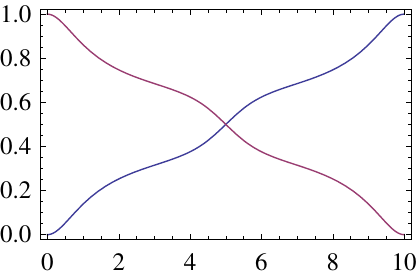}
\vspace*{-0.3cm}
\caption{Population change: $x^2_1 + x^2_2$ increasing, $x^2_3 + x^2_4$ decreasing.}
\end{center}
\end{figure}

The behaviour of the state variables are given in Fig. 5 indicating that the desired final state, where
$$x_3(T)=x_4(T)=0$$ has been obtained. The simulations were carried out and the figures produced by using Mathematica 7 package \cite{wolf7}.

\section{CONCLUSIONS}

The parametrization idea for constructing open-loop controls for uncontrollable bilinear systems is applied here. We have also earlier studied parametrization of systems described by partial differential equations and pseudo-differential operator models, see \cite{C46}--\cite{C48}. Flatness-based ideas, originally developed by Michel Fliess and his co-workers \cite{fliess92}-\cite{fliess99} have been developed for open-loop control design. In some quantum control problems, where laser pulses are used for the control, the dynamics is so fast that, at least at the present level of the speed of possible computations, feedback control seems to be impossible to implement even if so-called homodyne detection principles can be applied to obtain closed-loop controls.

Here we studied a two-level population transfer problem. Without more advanced differential geometric considerations, which might be helpful in understanding quantum phenomena in general, we use the formulation found generally in the literature, to obtain our basic driftless system model of the form $\dot x = g(x)u$, where $g$ is linear in the state $x$.

Simulation study was required to confirm the quantum control approach chosen. Then depending on the choise of the alternatives {\it a)} or {\it b)} different state trajectories can be obtained resulting, however, the same final state of the system when the flatness-based control is applied. Our preliminary simulations were based of the alternative {\it b)}.

The basic technique applied here is useful also in multi-qubit systems and in controlling entanglement of, say, two or more qubits. Then tensor product formalism in the Euclidian framework is a feasible alternative in the system model design.

\section{ACKNOWLEDGMENTS}
This work was supported in part by the European Commission, in Marie 
Curie programme's 
Transfer of Knowledge 
project {\it \mbox{Parametrization 
in the Control of Dynamic Systems} \mbox{(PARAMCOSYS, MTKD-CT-2004-509223)}}, which is greatly acknowledged.

\end{document}